\documentclass[a4paper,11pt,reqno]{amsart}
\usepackage[latin1]{inputenc}
\usepackage[english]{babel}
\usepackage{amssymb,amsfonts, amsmath, color}
\usepackage[margin=2.7cm]{geometry}
\usepackage[dvipsnames]{xcolor}
\usepackage{graphicx, color, enumerate}
\usepackage[latin1]{inputenc}
\usepackage[active]{srcltx}
\usepackage{tikz}
\usepackage{pgf}
\usepackage{etex}
\usepackage{verbatim}
\usepackage{tikz-3dplot}
\usepackage{pgfkeys}
\usepackage{amsmath}
\usepackage{amsfonts}
\usepackage{amssymb}
\usepackage{amsthm}
\usepackage{algpseudocode}
\usepackage{float}
\usepackage{xcolor}
\usepackage{mathrsfs}
\usepackage{import}
\usepackage{geometry}
\usepackage{fancyhdr}
\usepackage{fp}
\usepackage[colorlinks, citecolor=blue, linkcolor=red]{hyperref}

\newtheorem{theorem}{Theorem}[section]
\newtheorem{proposition}[theorem]{Proposition}
\newtheorem{corollary}[theorem]{Corollary}
\newtheorem{lemma}[theorem]{Lemma}
\newtheorem{remark}[theorem]{Remark}
\newtheorem{definition}[theorem]{Definition}

\newcommand{\bcl}{\begin{center}}
\newcommand{\ecl}{\end{center}}
\newcommand{\brl}{\begin{right}}
\newcommand{\erl}{\end{right}}
\newcommand{\ben}{\begin{enumerate}}
\newcommand{\een}{\end{enumerate}}
\newcommand{\overliner}{\begin{array}}
\newcommand{\earr}{\end{array}}
\newcommand{\btab}{\begin{tabular}}
\newcommand{\etab}{\end{tabular}}
\newcommand{\bdoc}{\begin{document}}
\newcommand{\edoc}{\end{document}}
\newcommand{\beqy}{\begin{eqnarray}}
\newcommand{\eeqy}{\end{eqnarray}}

\newcommand{\beqi}{\begin{eqnarray*}}
\newcommand{\eeqi}{\end{eqnarray*}}
\newcommand{\bitem}{\begin{itemize}}

\newcommand{\eitem}{\end{itemize}}
\newcommand{\nln}{\newline}
\newcommand{\newt}{\newtheorem}

\newcommand{\pa}{\partial}
\newcommand{\re}{{I\!\!R}}
\newcommand{\Rn}{\R^N}
\newcommand{\xr}{x\in\R }
\newcommand{\x}{\times}
\newcommand{\dyle}{\displaystyle}
\newcommand{\ene}{{I\!\!N}}
\newcommand{\irn}{\int\limits_{\R^N}}
\newcommand{\io}{\int\limits_{\O}}
\newcommand{\meas}{{\rm meas\,}}
\newcommand{\dif}{\nabla_{xy}}
\newcommand{\sign}{{\rm sign}}
\newcommand{\map}{\longrightarrow }
\newcommand{\imp}{\Longrightarrow }
\renewcommand{\div}{\nabla\cdot }
\newcommand{\sen}{{\rm sen\,}}
\newcommand{\tg}{{\rm tg\,}}
\newcommand{\arcsen}{{\rm arcsen\,}}
\newcommand{\arctg}{{\rm arctg\,}}
\newcommand{\supp}{{\textsl supp\ }}
\newcommand{\ity}{\int_{-\iy}^{+\iy}}
\newcommand{\limit}{\lim\limits}
\newcommand{\limi}{\limit_{n\to\infty}}
\newcommand{\sumi}{\sum\limits_{n=1}^{\infty}}
\newcommand{\ulu}{\underline u}
\newcommand{\ulw}{\underline w}
\newcommand{\ulz}{\underline z}
\newcommand{\ulv}{\underline v}
\newcommand{\uls}{\underline s}
\newcommand{\olu}{\overline u}
\newcommand{\olv}{\overline v}
\newcommand{\ols}{\overline s}
\newcommand{\ob}{\overline\b}
\newcommand{\ovar}{\overline\var}
\newcommand{\wv}{\widetilde v}
\newcommand{\wu}{\widetilde u}
\newcommand{\ws}{\widetilde s}
\renewcommand{\a }{\alpha }
\renewcommand{\b }{\beta }
\newcommand{\g }{\gamma}
\newcommand{\G }{\Gamma }
\renewcommand{\d }{\delta }

\newcommand{\D }{\Delta }
\newcommand{\e }{\varepsilon }
\newcommand{\z }{\zeta }
\renewcommand{\l }{\lambda }
\renewcommand{\L }{\Lambda }
\newcommand{\m }{\mu }
\newcommand{\n }{\nabla }
\newcommand{\s }{\sigma }
\newcommand{\Sig }{\Sigma }
\renewcommand{\t }{\tau }
\newcommand{\var }{\varphi }
\renewcommand{\o }{\omega }
\renewcommand{\O }{\Omega }
\newcommand{\R}{{\mathbb{R}}}
\newcommand{\bC}{{\bf C}}
\newcommand{\bZ}{{\bf Z}}
\newcommand{\bN}{{\bf N}}
\newcommand{\bQ}{{\bf Q}}
\newcommand{\bK}{{\bf K}}
\newcommand{\bI}{{\bf I}}
\newcommand{\bv}{{\bf v}}
\newcommand{\bV}{{\bf V}}
\newcommand{\LL}{\mathcal{L}}
\newcommand{\N}{\mathbb{N}}
\DeclareMathOperator{\suppo}{supp} \DeclareMathOperator{\di}{div}

\newenvironment{Proof}{\Rmovelastskip\vskip12pt
plus 1pt \noindent\em\rm}{\hfill {\qed \hskip .2cm}}

\begin{document}

\title[]{Phragm\`en-Lindel\"of type theorems \\ for  elliptic equations \\ on infinite graphs}

\author{Stefano Biagi}

\address{\hbox{\parbox{5.7in}{\medskip \noindent{Stefano Biagi, \\Dipartimento di Matematica, \\Politecnico di Milano, \\Piazza Leonardo da Vinci 32, 20133, Milano, Italy \\ [3pt] \emph{E-mail address: }{\tt stefano.biagi@polimi.it}}}}}



\author{Fabio Punzo}

\address{\hbox{\parbox{5.7in}{\medskip \noindent{Fabio Punzo, \\Dipartimento di Matematica, \\Politecnico di Milano, \\Piazza Leonardo da Vinci 32, 20133, Milano, Italy \\ [3pt] \emph{E-mail address: }{\tt fabio.punzo@polimi.it}}}}}



\begin{abstract} We investigate the validity of the Phragm\`en-Lindel\"of principle for a class of elliptic equations with a potential, posed on infinite graphs. Consequently, we get uniqueness, in the class of solutions satisfying a suitable growth condition at infinity. We suppose that the {\it outer degree (or outer curvature)} of the graph is bounded from above, and we allow the potential to go to zero at infinity in a controlled way. Finally, we discuss the optimality of the conditions on the potential and on the outer degree on special graphs.
\end{abstract}

\maketitle

\section{Introduction} \label{sec0}

Let
$(G, \omega, \mu)$ be
 a fixed {\it infinite weighted} graph, with {\it edge-weight} $\omega$ and {\it node (or vertex) measure} $\mu$.
In this paper we study uniqueness of possibly {\it unbounded} solutions to elliptic equations with a potential of the following form:
\begin{equation}\label{e1}
\Delta u -V u= f \quad \text{ in }\; G,
\end{equation}
where the potential $V$ is a nonnegative function defined in $G$
and $\Delta$ denotes the Laplace operator on $G$.

Recently, the analysis of partial differential equations on infinite graphs have attracted great attention (see,
e.g., \cite{Grig2, KLW, Mu2}). Qualitative properties of solutions to elliptic equations have been addressed in
\cite{AS1, AS2, BMP, HJ, HK, LT, Mas, MP2, MoPuSo, PinaS} for elliptic equations, while in \cite{BCG, CGZ, EM, GT,
HMu, HuangKS, LW2, MoPuSo2, Mu, Wu} for parabolic equations. In particular, $\ell^p$-Liouville theorems for
equation
\begin{equation}\label{e1h}
\Delta u -V u= 0 \quad \text{ in }\; G,
\end{equation}
with $V\equiv 0$ have been established in \cite{HK} (see also \cite{AS2, Mas, HJ}, and \cite{HKLS} for regular
Dirichlet forms). Furthermore, Liouville type results for equation \eqref{e1h}, with a potential $V$, can also be
deduced from Schnol's theorem for $p=2$; a discrete version of such a
result can be found in \cite{KLW} (see also \cite{LT}). In \cite{MP2}, Liouville type theorems have been stated
for equation \eqref{e1h} in {\it weighted} $\ell^p$ spaces, supposing that $\inf_G V >0.$ Such hypothesis has been
relaxed in \cite{BMP}, where it is assumed that $V$ can decay at infinity with a certain rate, but only {\it
bounded} solutions are considered.

The aim of this paper is to obtain a uniqueness result for equation \eqref{e1}, for possible unbounded solutions,
allowing the potential $V$ to decay at infinity, even if in a controlled way.  More precisely, we obtain
a Phragm\`en-Lindel\"of principle for equation \eqref{e1h} (see Proposition \ref{teo1}, Theorems \ref{cor1},
\ref{teo1z}). Then a Liouville type theorem for equation \eqref{e1h} and uniqueness of solutions for equation
\eqref{e1} immediately follow (see Corollaries \ref{cor2}, \ref{cor3}, \ref{cor4}).

There is a huge literature regarding Phragm\`en-Lindel\"of theorems for linear elliptic equations
in $\mathbb R^n$. We limit ourselves to quote \cite{Landis, Miller, Odd, OR, PoPuTe, ProtterW, PuTe}.

First we state a general Phragm\`en-Lindel\"of principle (see Proposition \ref{teo1}),
in which the existence of a suitable subsolution is assumed. Thus it is implicit in its character. Then we show
that on a class of graphs, such a subsolution can be constructed. Hence we obtain an explicit
Phragm\`en-Lindel\"of theorem (see Theorem \ref{cor1}), supposing that the potantial $V$ can decay at infinity
with a certain rate. Moreover, we assume that on such graphs, a relevant quantity, called {\it outer degree} (or
{\it outer curvature}), is bounded from above. Then we see that the hypothesis on $V$ and on the outer degree are optimal on spherically
symmetric trees. On the other hand, on the integer lattice $\mathbb Z^n$ we can allow a faster decay of the
potential $V$ (see Theorem \ref{teo1z}), and we show that such condition is sharp (see Theorem \ref{teo2z}).

\smallskip

The paper is organized as follows. In Section \ref{mf} the mathematical background is introduced,
then in Section \ref{sec2} and main results are stated. In Section \ref{auxiliary} we recall or obtain some
preliminary results. The
general Phragm\`en-Lindel\"of principle is proved in Section \ref{apriori};
finally, in Section \ref{integer} the results posed on the integer lattice $\mathbb Z^n$ are shown.

\section{Mathematical framework}\label{mf}\setcounter{equation}{0}

\subsection{The graph setting}
Let $G$ be a countably infinite set, and
let $\mu:G\to (0,+\infty)$ be a measure on $G$
satisfying $\mu(\{x\}) <+\infty$ for every $x\in G$ (so that $(G,\mu)$ becomes a measure space). Furthermore, let
\begin{equation*}
\omega:G\times G\to [0,+\infty)
\end{equation*}
be a symmetric, with zero diagonal and finite sum function, i.e.
\begin{equation}\label{omega}
\begin{aligned}
&\text{(i)}\,\, \omega(x,y)=\omega(y,x)\quad &\text{for all}\,\,\, (x,y)\in G\times G;\\
&\text{(ii)}\,\, \omega(x,x)=0 \quad\quad\quad\,\, &\text{for all}\,\,\, x\in G;\\
&\text{(iii)}\,\, \displaystyle \sum_{y\in G} \omega(x,y)<\infty \quad &\text{for all}\,\,\, x\in G\,.
\end{aligned}
\end{equation}
Thus, we define  \textit{weighted graph} the triplet $(G,\omega,\mu)$, where $\omega$ and $\mu$ are the so called \textit{edge weight} and \textit{node measure}, respectively. Observe that assumption $(ii)$ corresponds to ask that $G$ has no loops.
\smallskip

\noindent Let $x,y$ be two points in $G$; we say that
\begin{itemize}
\item $x$ is {\it connected} to $y$ and we write $x\sim y$, whenever $\omega(x,y)>0$;
\item the couple $(x,y)$ is an {\it edge} of the graph and the vertices $x,y$ are called the {\it endpoints} of the edge whenever $x\sim y$;
\item a collection of vertices $ \{x_k\}_{k=0}^n\subset G$ is a {\it path} if $x_k\sim x_{k+1}$ for all $k=0, \ldots, n-1.$
\end{itemize}

\noindent We are now ready to list some properties that the weighted graph $(G,\omega,\mu)$ may satisfy.

\begin{definition}\label{def01}
We say that the weighted graph $(G,\omega,\mu)$ is
\begin{itemize}
\item[(i)] {\em locally finite} if each vertex $x\in G$ has only finitely many $y\in G$ such that $x\sim y$;
\item[(ii)] {\em connected} if, for any two distinct vertices $x,y\in G$ there exists a path joining $x$ to $y$;
\end{itemize}
\end{definition}
For any $x\in G$, we define
\begin{itemize}
\item the {\it degree} of $x$ as $$\operatorname{deg}(x):=\sum_{y\in G}\omega(x,y);$$
\item the {\it weighted degree} of $x$ as $$\operatorname{Deg}(x):=\frac{\operatorname{deg}(x)}{\mu(x)}.$$
\end{itemize}

Let $d:G\times G\to [0, +\infty)$ be a distance on $G.$
For any $x_0\in G$ and $r>0$ we define the ball $B_r(x_0)$ with respect to the metric $d$ as
\[B_r(x_0):=\{x\in G\,:\, d(x,x_0)<r\}\,.\]
Sometimes we shall use specific metrics.

\noindent In this paper, we always make the following assumption:
\begin{equation}\label{e7f}
(G, \omega, \mu) \text{ is a connected, locally finite, weighted graph}.
\end{equation}

\subsection{Difference and Laplace operators} Let $\mathfrak F$ denote the set of all functions $f: G\to \mathbb R$\,. For any $f\in \mathfrak F$ and for all $x,y\in G$, let us give the following
\begin{definition}\label{def1}
Let $(G, \omega,\mu)$ be a weighted graph. For any $f\in \mathfrak F$,
\begin{itemize}
\item the {\em difference operator} is
\begin{equation*}
\nabla_{xy} f:= f(y)-f(x)\,;
\end{equation*}
\item the {\em (weighted) Laplace operator} on $(G, \omega, \mu)$ is
\begin{equation}\label{e2f}
\Delta f(x):=\frac{1}{\mu(x)}\sum_{y\in G}[f(y)-f(x)]\omega(x,y)\quad \text{ for all }\, x\in G\,.
\end{equation}
\end{itemize}
\end{definition}
Clearly,
\[\Delta f(x)=\frac 1{\mu(x)}\sum_{y\in G}(\dif f)\omega(x,y)\quad \text{ for all } x\in G\,.\]



\subsection{Outer and inner degrees}
We introduce some basic definitions following \cite[Chapter 9]{KLW}. Let $\rho$ be the {\it combinatorial graph distance} on $G$, that is, for any $x, y\in G,$ $\rho(x,y)$ is the least number of edges in a path between the two vertices $x$ and $y$.

Let $\Omega\subset G$ be finite subset. Define the distance from any $x\in G$ to the subset $\Omega$
\[\rho(x, \Omega):=\min_{y\in \Omega} \rho(x, y)\quad \forall x\in G\,.\]
With an abuse of notation we write $\rho(x, y)$ to indicate the distance between the two points $x, y\in G$, and $\rho(x, \Omega)$ to denote the distance from the point $x\in G$ to the set $\Omega\subset G$.

For any $r\in \mathbb N_0$, let
\[S_r(\Omega):=\{x\in G\,: \rho(x, \Omega)=r\,\}.\]
For any $f\in \mathfrak F$, if
\[f(x)=f(y) \quad \text{ whenever } \rho(x, \Omega)=\rho(y, \Omega),\]
then we say that $f$ is {\it spherically symmetric w.r.t.} $\Omega$; moreover, with a slight abuse of notation
we write
\[f(x)=f(r) \quad \forall x\in S_r(\Omega)\,.\]
For any $x\in G$ with $r\equiv r(x):=\rho(x, \Omega)\geq 1$, let
\[\mathfrak D_+(x):=\frac1{\mu(x)}\sum_{y\in S_{r+1}(\Omega)}\omega(x,y), \quad \mathfrak D_-(x):=\frac1{\mu(x)}\sum_{y\in S_{r-1}(\Omega)}\omega(x,y)\,.\]
The function $\mathfrak D_+:G\to [0, +\infty)$ is called {\it outer degree} w.r.t. $\Omega$, whereas $\mathfrak D_-:G\to [0, +\infty)$ is called {\it inner degree} w.r.t. $\Omega.$
Sometimes the outer degree and the inner degree are also called {\it outer curvature} and, respectively, {\it inner curvature} (see \cite{AS1}).

\begin{remark} \label{rem:omegaxySr}
We explicitly highlight, for a future reference,
 the following simple property: \emph{given any $x\in S_r(\Omega)$
(for some $r\geq 1$) and any $y\in G$, by definition of $\rho$ we have}
\begin{equation} \label{eq:omegaSr}
 \omega(x,y) > 0\,\,\Rightarrow\,\,y\in S_{r-1}(\Omega)\cup S_{r+1}(\Omega).
\end{equation}
\end{remark}

\medskip

The weighted graph $(G, \mu, \omega)$, endowed with the combinatorial distance $\rho$, is said to be
{\em weakly spherically symmetric} with
respect to a finite subset $\Omega\subset G$, if the outer and inner degrees $\mathfrak D_{\pm}$ are spherically symmetric with respect to $\Omega$. Therefore, on a weakly symmetric graph,
\[\mathfrak D_{\pm}(x)=\mathfrak D_{\pm}(r)\quad \forall x\in S_r(\Omega)\,.\]

\section{Main Results} \label{sec2}\setcounter{equation}{0}

\noindent We have already stated in \eqref{e7f}
the main hypotheses on the weighted graph $(G,\omega,\mu)$.
In order to state our main results, we first fix the following definition.
\begin{definition}\label{defsol}
Let $\mathcal{O}\subseteq G$ be an arbitrary set.
\begin{itemize}
 \item[\textbf{1)}] Given any $f:\mathcal{O}\to \mathbb{R}$, we say
 that a function $u\in\mathfrak{F}$ is a subsolution
 \emph{[}resp.\,su\-per\-so\-lu\-tion\emph{]} of the equation
 $\mathcal{L}u = f$ in $\mathcal{O}$ if, \emph{for every $x\in\mathcal{O}$}, we have
\begin{equation}\label{sol}
\frac{1}{\mu(x)}\sum_{y\in G}\omega(x,y)\left[u(y)-u(x)\right]-V(x)u(x)
\geq\,[\text{resp.\,\,$\leq$}]\,\,f(x).
\end{equation}
Moreover, we say that $u$ is a solution of the equation $\LL u = f$ in $\mathcal{O}$
if it is both a subsolution and a supersolution of the same equation.
\vspace*{0.1cm}

\item[\textbf{2)}] Given any $f:\mathcal{O}\to\mathbb{R}$ and $g:G\setminus\mathcal{O}\to\mathbb{R}$,
we say
 that a function $u\in\mathfrak{F}$ is a subsolution
 \emph{[}resp.\,supersolution\emph{]} of the
 $\LL$\,-\,Dirichlet problem
 \begin{equation} \label{eq:pbDir}
  \begin{cases}
 \LL u = f & \text{in $\mathcal{O}$} \\
 u = g & \text{in $G\setminus\mathcal{O}$}
 \end{cases}
 \end{equation}
 if the following conditions hold:
 \begin{itemize}
  \item[i)] $u$ is a subsolution \emph{[}resp.\,supersolution\emph{]} of the equation $\LL u = f$
  in $\mathcal{O}$;
  \item[ii)] $u\leq g$ \emph{[}resp.\,$u\geq g$\emph{]} pointwise in $G\setminus\mathcal{O}$.
 \end{itemize}
Finally, we say that $u$ is a solution of problem \eqref{eq:pbDir}
$u$ is both  a subsolution and a supersolution of this problem.
\end{itemize}
\end{definition}

\subsection{Phragm\`en-Lindel\"of principle}
The first main result
of this paper is a general Phragm\`en-Lindel\"of type principle, which reads as follows.

\begin{proposition}\label{teo1}
Let assumption \eqref{e7f} be satisfied. Let $V\in \mathfrak F, V\geq 0, x_0\in G$. Suppose that there exists $Z\in \mathfrak F$ bounded from above such that
\begin{equation}\label{e1p}
\Delta Z \geq - V \quad \text{ in } G,
\end{equation}
\begin{equation}\label{e2p}
\lim_{d(x, x_0)\to +\infty} Z(x) = -\infty\,.
\end{equation}
Let $u$ be a subsolution of equation \eqref{e1h} fulfilling
\begin{equation}\label{e3p}
\limsup_{d(x,x_0)\to +\infty}\frac{u(x)}{|Z(x)|} \leq 0\,.
\end{equation}
Then
\[u\leq 0 \quad \text{ in } G\,.\]
\end{proposition}

In particular, on  graphs with outer degree bounded from above, we have the next

\begin{theorem}\label{cor1}
Let assumption \eqref{e7f} be satisfied. Let $\Omega\subset G$ be a finite subset. Suppose that, for some $k_0>0$,
\begin{equation}\label{e15p}
0\leq \mathfrak D_+(x) \leq k_0\quad \text{ for all }\, x\in G \text{ with } \rho(x, \Omega)\geq 1\,.
\end{equation}
Suppose that
\begin{equation}\label{e12f}
V\in \mathfrak F, \quad V(x)\geq c_0\, [1+\rho(x, \Omega)]^{-\alpha} \quad \textrm{for all}\;\; x\in G,
\end{equation}
for some $\alpha\in [0, 1], c_0>0.$
 Let $u$ be a subsolution of equation \eqref{e1h} fulfilling
\begin{equation}\label{e16p}
\limsup_{\rho(x,\Omega)\to +\infty}\frac{u(x)}{|\tilde Z(x)|} \leq 0\,,
\end{equation}
where
\begin{equation}\label{e40f}
  \tilde Z(x):=
\begin{cases}
- [\rho(x, \Omega)]^{1-\alpha} & \text{ if } \alpha\in [0, 1)\\
-\log \rho(x, \Omega) & \text{ if } \alpha =1
\end{cases}, \text{ whenever } \rho(x, \Omega)>1.
\end{equation}
Then
\[u(x) \leq 0 \quad \forall x\in G.\]
\end{theorem}

We can immediately infer the following uniqueness results.

\begin{corollary}\label{cor2}
Let assumption \eqref{e7f} be satisfied. Let $V\in \mathfrak F, V\geq 0, x_0\in G$; let $f\in \mathfrak F$. Suppose that there exists $Z\in \mathfrak F$ bounded from above such that \eqref{e1p} and \eqref{e2p} hold. Then there exists at most one solution to equation \eqref{e1} such that
\[\lim_{d(x, x_0)\to +\infty}\frac{u(x)}{Z(x)}=0\,.\]
\end{corollary}

\begin{corollary}\label{cor3}
Let assumption \eqref{e7f} be satisfied; let $f\in \mathfrak F$. Assume that \eqref{e15p} and \eqref{e12f}  are fulfilled. Then there exists at most one solution to equation \eqref{e1} such that
\[\lim_{d(x, \Omega)\to +\infty}\frac{u(x)}{\tilde Z(x)}=0\,,\]
where $\tilde Z$ is given by \eqref{e40f}.
\end{corollary}

\subsection{Optimality on spherically symmetric trees}
In this subsection we consider a special kind of weakly symmetric graphs, the so called \textit{spherically symmetric trees}, and we show that the results in Theorem \ref{cor1} and
Corollary \ref{cor3} are sharp. More precisely, we show that \emph{both}
the choice $\alpha\in[0,1]$ in assumption \eqref{e12f}
and the boundedness assumption \eqref{e15p} on
$\mathfrak{D}_+$ \emph{are optimal}, that is,
infinitely many bounded solutions exist whenever $\alpha>1$ or whenever $\mathfrak{D}_+$ is not bounded.
\vspace*{0.1cm}

Let $(G, \omega, \mu)$ be a weakly symmetric graph w.r.t. $\Omega=\{o\}$, for some
fixed point $o\in G$
(which is usually referred to as the \emph{root of $G$}).
Suppose that
\begin{itemize}
\item $\omega:G\times G\to \{0, 1\}$;
\vspace*{0.1cm}

\item $\omega|_{S_r(\Omega)\times S_r(\Omega)}=0$;
\vspace*{0.1cm}

\item $\mu(x)=1 \text{ for every } x\in G;$
\vspace*{0.1cm}

\item there exists $b:\mathbb N\to\mathbb{N}$, which is called the {\it branching function}, such that
\[\mathfrak D_+(x)=b(r), \quad \mathfrak D_-(x)=1\quad \text{for every $x\in S_r(\Omega)$ and
$r\in \mathbb N$}.\]
\end{itemize}
As regards the optimality
of the choice $\alpha\in[0,1]$
in \eqref{e12f}, from \cite[Corollary 7.4]{BMP} we immediately
have the following result.
\begin{corollary}\label{nonuni}
  Let  $(G, \omega, \mu)$ be a spherically symmetric tree as above, with
  \emph{constant branching function} $b(r) = b_0\geq 2$. Moreover, let $V\in\mathfrak{F},\,V> 0$
on $G$. Assume that,

\begin{equation}\label{e8k}
V(x) \leq C_0\, \rho^{-\alpha}(x,\Omega)\quad \text{ for any }\; x\in G\setminus B_{R_0}\,.
\end{equation}
for some $\alpha>1$ and some $R_0\geq 2$.

Then
there exist infinitely many bounded solutions $u$ of problem \eqref{e1h}. In particular, for any $\gamma\in \R$, $\gamma>0$, there exists a solution $u$ to equation \eqref{e1h} such that
$$
\lim_{x\to\infty} u(x)=\gamma.
$$
\end{corollary}
As regards the optimality of the bounded assumption
\eqref{e15p}, instead, by relying on a general result
established in \cite{BMP} (see, precisely, Proposition \ref{prop71} below)
we will establish in the next
Section \ref{apriori} the following theorem.
\begin{theorem} \label{thm:OptimalityDplus}
 Let  $(G, \omega, \mu)$ be a spherically symmetric tree, with
 \emph{branching function}
 $$\text{$b(r) = (r+1)^p$,\quad for some $p\geq 1$.}$$
 Moreover, let  $0\leq \alpha\leq 1$ and let
 $$V:G\to\R,\qquad V(x) = (1+\rho(x,\Omega))^{-\alpha}.$$
 If $\alpha+p>1$, then
there exist infinitely many bounded solutions $u$ of problem \eqref{e1h}. In par\-ti\-cular, for any $\gamma\in \R$, $\gamma>0$, there exists a solution $u$ to equation \eqref{e1h} such that
$$
\lim_{x\to\infty} u(x)=\gamma.
$$
\end{theorem}
\begin{remark} \label{rem:comparisonResult}
 Let  $(G, \omega, \mu)$ be a spherically symmetric tree as above, with branching function
 $b(r)$. Moreover, let $\alpha\in[0,1]$ be fixed, and let
 $$V:G\to\R,\qquad V(x) = (1+d(x,\Omega))^{-\alpha}.$$
 On account of Theorem \ref{cor1} we see that, \emph{if $b(r)$ is bounded},
 then there exists at most one bounded solution
 to equation \eqref{e1h}; on the other hand, if
 $b(r) = (1+r)^p$ (for some $p\geq 2$),
 it follows from Theorem \ref{thm:OptimalityDplus} that there exists
 \emph{infinitely many bounded solutions} of \eqref{e1h}.
\end{remark}

\subsection{Further results on $\mathbb Z^n$}
Consider the $n-$dimensional {\it integer lattice graph}. Then $G=\mathbb Z^n$; furthermore, $x\sim y$ if and only if there exists $k\in \{1, \ldots, n\}$ such that $x_k=y_k\pm 1$ and $x_i=y_i$ for $i\neq k$. The edge weight is
\[\omega: \mathbb Z^n\times \mathbb Z^n\to [0, +\infty)\]
\[\omega(x,y)=\begin{cases}
1 & \text{ if } y\sim x\\
0 & \text{ if } y\not\sim x\,,
\end{cases}\]
whereas the node measure is
\[\mu(x)=\sum_{y\in \mathbb Z^n} \omega(x,y)=2n\,.\]

We equip the graph $(\mathbb Z^n, \omega, \mu)$ with the euclidean distance
\[|x-y|=\left(\sum_{k=1}^n |x_k-y_k|^2 \right)^{\frac 12}\quad (x,y\in \mathbb Z^n)\,.\]

\begin{remark}
We remark that $\mathbb Z^n$ with the euclidean distance is not a weakly symmetric graph, since in the definition of weakly symmetric graph only the combinatorial graph distance is considered. Furthermore, it is easily seen that $\mathbb Z^n$ endowed with the combinatorial metric is not a weakly symmetric graph.
\end{remark}

In this case the condition on $\alpha$ made in \eqref{e12f} is not optimal.
In fact, the critical value is now $\alpha=2$, and not more $\alpha=1$.

More precisely, we can prove the next results.
\begin{theorem}\label{teo1z}
Let $G=\mathbb Z^n$. Suppose that, for some $\alpha\in [0, 2]$ and $c_0>0$,
\begin{equation}\label{e12ff}
V\in \mathfrak F, \quad V(x)\geq c_0\, (1+|x|)^{-\alpha} \quad \textrm{for all}\;\; x\in G.
\end{equation}
Let $u$ be a subsolution of equation \eqref{e1h} fulfilling
\eqref{e3p},
with $x_0=0$, $d(x,y)=|x-y|$ and
\begin{equation}\label{e40ff}
 Z(x):=
\begin{cases}
-|x|^{2-\alpha} & \text{ if } \alpha\in [0, 2)\\
-\log |x| & \text{ if } \alpha =2
\end{cases}, \text{ whenever } |x|>2.
\end{equation}
Then
\[u(x) \leq 0 \quad \forall x\in G.\]
\end{theorem}

A direct consequence of Theorem \ref{teo1} is the following uniqueness result.
\begin{corollary}\label{cor4}
Let $G=\mathbb Z^n, f\in \mathfrak F$. Assume that \eqref{e12ff} holds. Then there exists at most one solution to equation \eqref{e1} such that
\[\lim_{|x|\to +\infty}\frac{u(x)}{Z(x)}=0\,,\]
where $Z$ is given by \eqref{e40ff}.
\end{corollary}

The value $\alpha=2$ in Corollary \ref{cor4} is optimal. This follows from the next nonuniqueness result.

\begin{theorem}\label{teo2z}
Let $G=\mathbb Z^n, n\geq 3$.  Assume that
\begin{equation}\label{e12fb}
V\in \mathfrak F, \quad 0<V(x)\leq c_0\, (1+|x|)^{-\alpha} \quad \textrm{for all}\;\; x\in G,
\end{equation}
for some $\alpha>2$. Then there exist infinitely many bounded solutions to equation \eqref{e1h}. More precisely, for any $\gamma\in \mathbb R$ there exists a bounded solution $u_\gamma$ such that
\[\lim_{|x|\to +\infty} u_\gamma(x)=\gamma\,.\]
\end{theorem}

\section{Auxiliary Results}\label{auxiliary}\setcounter{equation}{0}
In this section we collect some well-known facts
which will be exploited in the proofs of our main results.
To begin with, we prove the following simple lemma.
\begin{lemma}
Let assumption \eqref{e7f} be satisfied. Let $\Omega\subset G$ be a finite set and let $f\in \mathfrak F$ be a spherically symmetric function with respect to $\Omega$.
Then
\begin{equation}\label{e13p}
\Delta f(x)=\mathfrak D_+(x)[f(r+1)-f(r)]+\mathfrak D_-(x)[f(r-1)-f(r)]
\end{equation}
for any $x\in G$ with $r\equiv r(x):=\rho(x, \Omega)\geq 1.$
\end{lemma}
\begin{proof}
Let $x\in G$ with $r\equiv r(x):=\rho(x, \Omega)\geq 1.$
From \eqref{e2f} and Remark \ref{rem:omegaxySr} we have
\begin{equation}\label{e71f}
\begin{aligned}
\Delta f(x)&=\Delta f(r(x))=\frac 1{\mu(x)}\sum_{y\in G}[f(y)-f(x)]\omega(x,y)\\
&=\frac 1{\mu(x)}\sum_{y\in S_{r+1}(\Omega)\cup S_{r-1}(\Omega)}[f(y)-f(x)]\omega(x,y)\\
&= \frac 1{\mu(x)}\Big\{ \sum_{y\in S_{r+1}(\Omega)}[f(y)-f(x)]\omega(x,y)\\
& \qquad\qquad+ \sum_{y\in S_{r-1}(\Omega)}[f(y)-f(x)]\omega(x,y)\Big\}\\
&= \frac{1}{\mu(x)}\sum_{y\in S_{r+1}(\Omega)}[f(r+1)-f(r)]\omega(x,y) \\
&\qquad\qquad+  \frac{1}{\mu(x)}\sum_{y\in S_{r-1}(\Omega)}[f(r-1)-f(r)]\omega(x,y)\\
&= [f(r+1)-f(r)]\frac{1}{\mu(x)}\sum_{y\in S_{r+1}(\Omega)}\omega(x,y) \\
&\qquad\qquad+ [f(r-1)-f(r)]\frac{1}{\mu(x)}\sum_{y\in S_{r-1}(\Omega)}\omega(x,y) \\
&=\mathfrak D_+(x)[f(r+1)-f(r)]+ \mathfrak D_-(x)[f(r-1)-f(r)]\,.
\end{aligned}
\end{equation}
This ends the proof.
\end{proof}
We then proceed by reviewing some results proved in \cite{BMP}. Throughout the rest
of the paper, to simplify the notation we denote by $\mathcal{L}$ the operator
$$\mathcal{L} := \Delta-V(x).$$
First of all, we state the following \emph{Weak Maximum Principle} (see, e.g., \cite[Lemma 3.3]{BMP}).
\begin{lemma} \label{lem:WMP} Let assumption \eqref{e7f} be in force.
 Let $\Omega\subseteq G$ be a finite set, and let $u\in\mathfrak{F}$ be
 a supersolution of the \emph{homogeneous Dirichlet problem}
 \begin{equation} \label{eq:PBWMP}
  \begin{cases}
   \LL u = 0 & \text{in $\Omega$} \\
   u = 0 & \text{in $G\setminus\Omega$}
  \end{cases}
 \end{equation}
 \emph{(}in the sense of Definition \ref{defsol}\emph{)}.
 Then $$u\geq 0 \quad \text{ in } \Omega.$$
\end{lemma}
Then, we recall a general non-uniqueness criterion (see \cite[Proposition 7.1]{BMP}).
\begin{proposition}\label{prop71}
Let assumption \eqref{e7f} be in force. Moreover, let $V\in\mathfrak{F}$, $V>0, x_0\in G $ and $\hat R>0$. Assume
that there exists a supersolution of the equation
\begin{equation*}
 \Delta h=-V \quad \text{in}\,\,\,\mathcal{O} = G\setminus B_{\hat R}(x_0),
\end{equation*}
further satisfying the following properties
\begin{itemize}
\item[i)]
$h>0$ in $G\setminus B_{\hat R}(x_0)$;
\vspace*{0.1cm}

\item[ii)] $\lim_{x\to\infty}h(x)=0.$
\end{itemize}
Then there exist infinitely many bounded solutions $u$ of problem \eqref{e1h}. In particular, for any $\gamma\in \R$, $\gamma>0$, there exists a solution $u$ to problem \eqref{e1h} such that
$$
\lim_{d(x, x_0)\to\infty} u(x)=\gamma.
$$
\end{proposition}

\section{Proof of Proposition \ref{teo1} and Theorem \ref{cor1}}\label{apriori}\setcounter{equation}{0}
\begin{lemma}\label{lemma1}
Let $Z$ be a subsolution of equation
\begin{equation}\label{e5p}
\Delta Z = - V \quad \text{ in } G\,,
\end{equation}
such that, for some $H\in \mathbb R$,
\begin{equation}\label{e6p}
Z\leq H \quad \text{ in } G\,.
\end{equation}
Then
\begin{equation}\label{e7p}
\bar Z:= Z- H -1
\end{equation}
fulfills
\begin{equation}\label{e9p}
\bar Z\leq -1 \quad \text{ in } G,
\end{equation}
\begin{equation}\label{e8p}
\Delta \bar Z - V \bar Z \geq 0 \quad \text{ in } G\,.
\end{equation}
\end{lemma}
\begin{proof}[Proof of Lemma \ref{lemma1}]
Obviously, \eqref{e7p} is an immediate consequence of \eqref{e6p}. Furthermore, from \eqref{e5p} and \eqref{e9p} we get
\begin{equation*}
\begin{aligned}
\Delta \bar Z - V \bar Z &= \Delta Z - V\bar Z \\ & \geq - V - V \bar Z = V(-1-\bar Z)\\ & \geq V (-1+1)=0 \quad \text{ in } G\,,
\end{aligned}
\end{equation*}
hence \eqref{e8p} follows.
\end{proof}

\begin{proof}[Proof of Proposition \ref{teo1}] In view of the hypotheses, we can apply Lemma \ref{lemma1} to get the existence of $\bar Z$ fulfilling \eqref{e9p} and \eqref{e8p}.
From \eqref{e2p}, \eqref{e3p} and \eqref{e9p} we can infer that
\[\liminf_{d(x,x_0)\to +\infty} \frac{u(x)}{\bar Z(x)}\geq 0\,.\]
Thus there exists a sequence $\{R_n\}\subset (0, +\infty)$ with $\displaystyle R_n \mathop{\longrightarrow}_{n\to +\infty} +\infty$ such that
\[\lim_{n\to +\infty} \inf_{G\setminus B_{R_n}(x_0)}\frac{u(x)}{\bar Z(x)}\geq 0\,.\]
This means that for any $\alpha>0$ there exists $\bar n=\bar n(\alpha)\in \mathbb N$ such that, for all $n>\bar n$,
\[
\frac{u}{\bar Z}>-\alpha \quad \text{ in } G\setminus B_{R_n}(x_0),
\]
thus
\begin{equation}\label{e10p}
u<-\alpha \bar Z \quad \text{ in } G\setminus B_{R_n}(x_0)\,.
\end{equation}

\medskip

For any $\alpha>0$ define
\[\mathcal Z_\alpha:= -\alpha \bar Z\,.\]
By virtue of \eqref{e9p},
\begin{equation}\label{e11p}
\mathcal Z_\alpha \geq \alpha \quad \text{ in } G\,.
\end{equation}
From \eqref{e8p} and \eqref{e11p} it follows that for any $\alpha>0, n\in \mathbb N$, $\mathcal Z_\alpha$ is a supersolution of problem
\begin{equation}\label{e12p}
\begin{cases}
\LL u = 0 & \text{ in } B_{R_n}(x_0)\\
u = \mathcal Z_\alpha & \text{ in } G\setminus B_{R_n}(x_0)\,.
\end{cases}
\end{equation}
On the other hand, thanks to \eqref{e10p}, for any $\alpha>0$ and for any $n>\bar n$ (with $\bar n\in \mathbb N$  as above), $u$ is a subsolution of problem \eqref{e12p}.
By Lemma \ref{lem:WMP},
\[u\leq \mathcal Z_\alpha \quad \text{ in } B_{R_n}(x_0)\,.\]
Letting $\alpha\to 0^+$, we deduce that
\[u\leq 0 \quad \text{ in } G\,.\]
\end{proof}

\begin{proof}[Proof of Theorem \ref{cor1}]
Let
\[r\equiv r(x):= \rho(x, \Omega)\quad \forall\, x\in G\,. \]
To begin with, we consider the case where \eqref{e12f} holds with $\alpha\in [0, 1)\,.$

Define
\[Z(x)\equiv Z(r):=- M r^\beta - 1 \quad \forall x\in G\,,\]
where $M>0$ will be properly determined in the sequel, whereas $\beta:=1-\alpha$.

In view of \eqref{e13p}, for all $x\in G, r(x)\geq 1$,
\begin{equation*}
\begin{aligned}
\Delta Z(x)&=M \left\{-\mathfrak D_+(x)[(r+1)^\beta-r^\beta] + \mathfrak D_-(x)[r^\beta- (r-1)^\beta]\right\}\\
&\geq M \left\{-\mathfrak D_+(x)[(r+1)^\beta-r^\beta] \right\} = M\left\{-\mathfrak D_+(x) \beta \xi^{\beta-1} \right\},
\end{aligned}
\end{equation*}
for some $\xi\in (r, r+1)$. Therefore, since $\beta<1$,
\begin{equation}\label{e17p}
\Delta Z(x) \geq - \frac{M \beta k_0}{r^{1-\beta}} = - \frac{M \beta k_0}{r^{\alpha}} \quad \forall x\in G, r(x)\geq 1\,.
\end{equation}
Combining together \eqref{e12f} and \eqref{e17p} we can infer that
\begin{equation}\label{e19p}
\frac 1{V(x)}\Delta Z(x) \geq - \frac 1{c_0}\frac{(1+r)^\alpha}{r^{\alpha}}M \beta k_0\geq -1 \quad \forall x\in G, r(x)\geq 1\,,
\end{equation}
provided that $M>0$ is small enough. On the other hand, since $\Omega$ is a finite subset of $G$ (and since $V>0$), by shrinking $M>0$ if needed, we also have that
\begin{equation}\label{e20p}
\Delta Z(x) = -M\Delta(r^\beta) \geq - V(x)\quad\forall\,\,x\in\Omega.
\end{equation}
By virtue of \eqref{e19p} and \eqref{e20p}, $Z$ fulfills \eqref{e1p}.

Now, let $d_0$ be the diameter of the finite set $\Omega$, let $x_0\in \Omega$. Select any $x\in G$ with $\rho(x, x_0)\geq 2d_0$. For all $y\in \Omega$, by triangular inequality,
\[\rho(x,y)\geq \rho(x, x_0) - \rho(y, x_0)\geq \rho(x, x_0) - d_0.\]
Hence
\[\rho(x, \Omega)=\min_{y\in \Omega} \rho(x,y)\geq \rho(x, x_0) - d_0,\]
thus
\begin{equation}\label{e18p}
\rho(x, x_0)- d_0\leq \rho(x, \Omega)\leq \rho(x, x_0)\,.
\end{equation}
In view of \eqref{e18p}, $Z$ satisfies also \eqref{e2p}. By \eqref{e18p} again, since by assumption $u$ satisfies \eqref{e16p}, we can infer that also \eqref{e3p} holds.
Therefore we are in position to apply Theorem \ref{teo1}, with $d=\rho$, to get the thesis, when $\alpha\in [0, 1)$.

\medskip

It remains to take into account the case where $\alpha=1$.
Define
\[Z(x)\equiv Z(r):=- M \log(2+r)  \quad \forall x\in G\,,\]
where $M>0$ will be properly determined in the sequel.

In view of \eqref{e13p}, for all $x\in G, r(x)\geq 1$,
\begin{equation*}
\begin{aligned}
\Delta Z(x)&=M \left\{-\mathfrak D_+(x)[-\log(2+r)+\log(3+r)] + \mathfrak D_-(x)[\log(2+r)-\log(1+r) ]\right\}\\
&\geq M \left\{-\mathfrak D_+(x)[-\log(2+r)+\log(3+r)] \right\} = -\frac{M}{\xi} \mathfrak D_+(x)   ,
\end{aligned}
\end{equation*}
for some $\xi\in (2+r, 3+r)$. Therefore,
\begin{equation}\label{e20pBIS}
\Delta Z(x) \geq - \frac{M  k_0}{2+r} \quad \forall r\in \mathbb N_0\,.
\end{equation}
Combining together \eqref{e12f} and \eqref{e20pBIS} we can infer that
\begin{equation}\label{e21p}
\frac 1{V(x)}\Delta Z(x) \geq - \frac 1{c_0}\frac{1+r}{2+r}M  k_0\geq -1 \quad \forall x\in G, r(x)\geq 1,
\end{equation}
provided that $M>0$ is small enough. On the other hand, by proceeding exactly as in the previous case (and possibly shrinking the constant $M$), we also have that
\begin{equation}\label{e22p}
\Delta Z(x)=-M\Delta(\log(2+r))\geq - V(x)\quad\forall\,\,x\in\Omega.
\end{equation}
By virtue of \eqref{e21p} and \eqref{e22p}, $Z$ fulfills \eqref{e1p}. Similarly to the previous case, we get the thesis by means of Theorem \ref{teo1}, with $d=\rho$.
This completes the proof.
\end{proof}
We conclude the section by proving Theorem \ref{thm:OptimalityDplus}.
\begin{proof}[Proof of Theorem \ref{thm:OptimalityDplus}]
 Let $\beta > 0$ be arbitrarily fixed, and let
 $$h:G\to\R, \qquad h(x) = (1+d(x,\Omega))^{-\beta}.$$
 We claim that it is possible to choose $\beta$ in such a way that
 \begin{equation} \label{eq:hsupersolOptimal}
 \frac{1}{V(x)}\Delta h(x)\leq -1\quad\text{for all $x\in G$ with $r = d(x,\Omega)\geq R_0$,}
 \end{equation}
 provided that $R_0>0$ is sufficiently large
 (recall that $\Omega = \{o\}$, where $o$ is the \emph{root of $G$}). Once this claim
 has been established,
 the thesis of the theorem (i.e., the existence
 of infinitely many solutions to problem \eqref{e1h}) is direct
 consequence of Proposition \ref{prop71} (note that $h$
 satisfies properties i)\,-\,ii) in the statement of this proposition).
 Hence, we turn to prove \eqref{eq:hsupersolOptimal}.
 \vspace*{0.1cm}

We first observe that, since we are assuming
 $\mathfrak{D}_+(x) = b(r) = (r+1)^p,\,
 \mathfrak{D}_-(x)\equiv 1$
 (for every fixed $x\in S_r(\Omega)$ and every $r\geq 1$),
a direct application of Lemma \ref{lem1z} gives
\begin{align*}
 \Delta h(x) & = (1+r)^{p}[(2+r)^{-\beta}-(1+r)^{-\beta}]
 + [r^{-\beta}-(1+r)^{-\beta}] \\
 & (\text{by the Mean Value theorem, for some $r+1<\xi<r+2$ and $r<\eta<r+1$})\\
 & = -\beta(1+r)^p\xi^{-\beta-1}+ \beta \eta^{-\beta-1} \\
 & (\text{since $\xi < r+2\leq 2r$, as $r\geq 2$}) \\
 & \leq -2^{-\beta-1}\beta r^{p-\beta-1}+\beta r^{-\beta-1}.
\end{align*}
From this, we obtain
\begin{equation} \label{eq:StimaFinaleOptimal}
\begin{split}
 \frac{1}{V(x)}\Delta h(x) & = -2^{-\beta-1}\beta (1+r)^{\alpha}r^{p-\beta-1}
 + \beta(1+r)^{\alpha}r^{-\beta-1} \\
 & \leq -2^{-\beta-1}\beta r^{p+\alpha-\beta-1}+ 2^\alpha \beta r^{\alpha-\beta-1},
 \end{split}
 \end{equation}
 and this estimate holds for every $x\in G$ with $r = d(x,\Omega)\geq 2$.
 In view of \eqref{eq:StimaFinaleOptimal},
 if we choose $\beta>0$ in such a way that
 $0<\beta<\alpha+p-1$ (recall that $\alpha+p> 1$), we conclude that
 $$\frac{1}{V(x)}\Delta h(x)\leq -2^{-\beta-1}\beta r^{p+\alpha-\beta-1}+ 2^\alpha\beta r^{\alpha-\beta-1}
 \leq -1,$$
 provided that $r = d(x,\Omega)$ is sufficiently large. This ends the proof.
\end{proof}

\section{Further results on $\mathbb Z^n$: proofs}\label{integer} \setcounter{equation}{0}
\subsection{Construction of sub/supersolutions on $\mathbb Z^n$}

\begin{lemma}\label{lem1z}
Let the assumptions of Theorem \ref{teo1z} hold with $0\leq \alpha<2$. Then, for some $\beta>0, K>0$,
\[Z(x):=- K |x|^{2\beta} - 1 \quad (x\in \mathbb Z^n)\]
fulfills
\begin{equation}\label{e60f}
\frac1{V(x)} \Delta Z(x) \geq -1 \quad \forall x\in \mathbb Z^n\,.
\end{equation}
\end{lemma}
\begin{proof}[Proof of Lemma \ref{lem1z}] Let
\[\varphi(t):=- K t^\beta -1 \quad \forall t\geq 0\,,\]
where $K>0, \beta>0.$
Obviously,
\[\varphi'(t)=-\beta K t^{\beta -1},\]
\[\varphi''(t)=\beta(1-\beta)K t^{\beta-2},\]
\[\varphi'''(t)=\beta(1-\beta)(\beta-2)K t^{\beta-3}\,,\]
\begin{equation}\label{e50f}
\varphi(t)=\varphi(r)+\varphi'(r)(t-r) +\frac{\varphi''(r)}{2}(t-r)^2+ \frac{\varphi'''(\eta)}{6}(t-r)^3
\end{equation}
for any $t, r\in (0, +\infty),$ for some $\eta$ between $t$ and $r$.
\vspace*{0.1cm}

By virtue of \eqref{e50f} with $t=|y|^2,\,r=|x|^2$ ($x, y\in \mathbb Z^n, |x|\geq 1$),
\begin{equation}\label{e51f}
\begin{aligned}
\Delta Z(x)&= \frac 1{\mu(x)}\sum_{x\in \mathbb Z^n}[Z(y)- Z(x)]\omega(x,y) \\ &
= \frac K{2n} \sum_{x\in \mathbb Z^n}\Big\{\varphi'(|x|^2)[|y|^2- |x|^2]+\frac{\varphi''(|x|^2)}{2}[|y|^2-|x|^2]^2\\
&\qquad\qquad+\frac{\varphi'''(\eta)}{6}[|y|^2-|x|^2]^3 \Big\}\omega(x,y)\\
&= \frac K{2n} \sum_{y\in \mathbb Z^n}\Big\{-\beta |x|^{2\beta-2}(|y|^2-|x|^2)+\beta(1-\beta)|x|^{2\beta-4}(|y|^2-|x|^2)^2\\
& \qquad\qquad+ \frac{\beta}{6}(1-\beta)(\beta-2)\eta^{\beta-3}(|y|^2-|x|^2)^3  \Big\}\omega(x,y),
\end{aligned}
\end{equation}
for some $\eta$ fulfilling
\begin{equation}\label{e54f}
\min\{|x|^2, |y|^2\}\leq \eta\leq \max\{|x|^2, |y|^2\}\,.
\end{equation}
We then turn to estimate from below the right-hand side of \eqref{e51f}. In what
follows, to ease the readability we denote by $c$ any positive constant, possibly different from line
to line but independent of the fixed $x,y\in G$ with $\omega(x,y) > 0$.
\vspace*{0.1cm}

Let $x\in \mathbb Z^n$ be fixed. Define
\[\mathcal O_x:=\{y\in \mathbb Z^n\,:\, \omega(x,y)>0\}.\]
It is easy to see that $\mathcal O_x$ contains $n$ points $y\in \mathbb Z^n$ of the form
$y=(x_1, \ldots, x_{k}- 1, \ldots, x_n)$ for any $k\in \{1, \ldots, n\}$, and $n$ points $y\in \mathbb Z^n$ of the form
$y=(x_1, \ldots, x_{k}- 1, \ldots, x_n)$ for any $k\in \{1, \ldots, n\}$.
Therefore, there are $n$ points $y\in \mathbb Z^n$ for which
\[|y|^2-|x|^2=|x|^2+ 2x_k +1-|x|^2= 2x_k+1\,,\]
and $n$ points $y\in \mathbb Z^n$ for which
\[|y|^2-|x|^2=|x|^2- 2x_k +1-|x|^2= - 2x_k+1\,,\]
Consequently,
\begin{equation}\label{e52f}
\sum_{y\in \mathbb Z^n}(|y|^2-|x|^2)\omega(x,y)=2n\,.
\end{equation}
Furthermore, we have
\begin{equation}\label{e53f}
\left|[|y|^2-|x|^2]^3 \right|\leq c(|x|^3+1) \quad \forall x,y\in \mathbb Z \text{ with } \omega(x,y)>0\,.
\end{equation}
Let now $x, y\in \mathbb Z^n$ for which $\omega(x,y)>0$, and let $\eta$ be as in \eqref{e54f}.
If $|x|\geq 1$, then $|y|\geq |x|-1$; as a consequence, we can write
\begin{equation}\label{e55f}
\eta\geq (|x|-1)^2\geq c(|x|^2-1).
\end{equation}
Let
\[\mathfrak R(x):=\frac{\beta}{6}(1-\beta)(\beta-2)\eta^{\beta-3}(|y|^2-|x|^2)^3.\]
Due to \eqref{e53f} and \eqref{e55f}, for every $x, y\in \mathbb Z^n$ for which $\omega(x,y)>0$  we have
\begin{equation}\label{e57f}
\begin{aligned}
\mathfrak R(x)&\geq -\frac{c\beta(1-\beta)(\beta-2)}{6}(|x|^2-1)^{\beta -3}(1+|x|^3)\\
&\geq -c(1+|x|)^{2\beta-3}\,.
\end{aligned}
\end{equation}
Hence, by virtue of \eqref{e57f} and \eqref{e52f}, choosing $0<\beta<1$, \eqref{e51f} yields
\begin{equation}\label{e56f}
\begin{aligned}
\Delta Z(x)& \geq  \frac K{2n} \sum_{y\in \mathbb Z^n}
 \left\{ - \beta  |x|^{2\beta-2}(|y|^2-|x|^2) +\mathfrak R(x)\right\}\omega(x,y) \\
 & \geq - K \left\{\beta |x|^{2\beta-2}+ c(1+|x|)^{2\beta-3}\right\}\omega(x,y) \\
& \geq - c K |x|^{2\beta-2}  \quad \forall x\in \mathbb Z^n, |x|\geq 1\,.
\end{aligned}
\end{equation}
Observe that in view of \eqref{e12ff}, for some $\alpha\in [0, 2)$,
\[\frac 1{V(x)} \leq \frac 1{c_0}(1+|x|)^{\alpha}\quad \forall x\in \mathbb Z^n\,.\]
Therefore, from \eqref{e56f} we can infer that
\begin{equation}\label{e57ff}
\frac{1}{V(x)}\Delta Z(x)\geq - c K |x|^{2\beta-2 +\alpha}\geq -1\quad \forall x\in \mathbb Z^n, |x|\geq 1,
\end{equation}
provided that $K\leq \frac1{c}$ and $0<\beta<\frac{2-\alpha}{2}.$
On the other hand, we also have
\begin{equation}\label{e58f}
\frac 1{V(0)}\Delta Z(0) = -K\frac{\Delta(|x|^{2\beta})(0)}{V(0)}\geq -1,
\end{equation}
by possibly shrinking $K$. Gathering \eqref{e57ff} and \eqref{e58f}, we get that $Z$ satisfies \eqref{e60f}.
\end{proof}
\begin{lemma}\label{lem2z}
Let the assumptions of Theorem \ref{teo1z} hold with $0\leq \alpha\leq 2$. Then, for some $K>0$,
the function
\[Z(x):=- K \log\big(|x|^2 +2\big) \quad (x\in \mathbb Z^n)\]
fulfills \eqref{e60f}.
\end{lemma}

\begin{proof}[Proof of Lemma \ref{lem2z}] Let
\[\psi(t):=- K \log(t+2) \quad \forall t\geq 0\,,\]
where $K>0.$
Obviously,
\[\psi'(t)=-\frac K{t+2},\]
\[\psi''(t)=\frac{K}{(t+2)^2},\]
\[\psi'''(t)=-\frac{2 K}{(t+2)^3}\,,\]
\begin{equation}\label{e50fb}
\psi(t)=\varphi(r)+\psi'(r)(t-r) +\frac{\psi''(r)}{2}(t-r)^2+ \frac{\psi'''(\eta)}{6}(t-r)^3
\end{equation}
for any $t, r\in (0, +\infty),$ for some $\eta$ between $t$ and $r$.
\vspace*{0.1cm}

By virtue of \eqref{e50fb} with $t=|y|^2, r=|x|^2\quad (x, y\in \mathbb Z^n, |x|\geq 1)$,
\begin{equation}\label{e51fb}
\begin{aligned}
\Delta Z(x)&= \frac 1{\mu(x)}\sum_{x\in \mathbb Z^n}[Z(y)- Z(x)]\omega(x,y) \\ &
= \frac K{2n} \sum_{x\in \mathbb Z^n}\Big\{\psi'(|x|^2)[|y|^2- |x|^2]+
 \frac{\psi''(|x|^2)}{2}[|y|^2-|x|^2]^2\\
 &\qquad\qquad+\frac{\psi'''(\eta)}{6}[|y|^2-|x|^2]^3 \Big\}\omega(x,y)\\
&=\frac K{2n} \sum_{y\in \mathbb Z^n}\Big\{-\frac{1}{|x|^2+2}(|y|^2-|x|^2)+\frac 1{2(|x|^2+2)^2}(|y|^2-|x|^2)^2\\ & \qquad\qquad- \frac 1{3(\eta+2)^3}(|y|^2-|x|^2)^3  \Big\}\omega(x,y),
\end{aligned}
\end{equation}
for some $\eta$ fulfilling \eqref{e54f}.
As in the proof of Lemma \ref{lem1z}
we then proceed by estimating from below the right-hand side of \eqref{e51fb}.
In what follows, we simply denote by $c$ any positive con\-stant, possibly different
from line to line but independent of $x,y$.

Let
\[\mathfrak T(x):=-\frac 1{3(\eta+2)^3}(|y|^2-|x|^2)^3.\]
Due to \eqref{e53f} and \eqref{e55f}, for every $x, y\in \mathbb Z^n$ for which $\omega(x,y)>0$,  we have
\begin{equation}\label{e57fb}
\begin{aligned}
\mathfrak T(x)&\geq -\frac{c}{3}(|x|^2+1)^{-3}(1+|x|^3)\\
&\geq - c(1+|x|)^{-3}\,.
\end{aligned}
\end{equation}
Hence, by virtue of \eqref{e57fb} and \eqref{e52f}, \eqref{e51fb} yields
\begin{equation}\label{e56fb}
\begin{aligned}
\Delta Z(x)& \geq  \frac K{2n} \sum_{y\in \mathbb Z^n}\left\{ -\frac{1}{|x|^2+2}(|y|^2-|x|^2) +\mathfrak T(x)\right\}\omega(x,y)\\ & \geq -  K \left\{\frac 1{|x|^2}+ \frac{c}{(1+|x|)^{3}}\right\}
\geq - c K |x|^{-2}  \quad \forall\,\,x\in \mathbb Z^n, |x|\geq 1\,.
\end{aligned}
\end{equation}
Observe that in view of \eqref{e12ff}, for some $\alpha\in [0, 2]$,
\[\frac 1{V(x)} \leq \frac 1{c_0}(1+|x|)^{\alpha}\quad \forall\,\,x\in \mathbb Z^n\,.\]
Therefore, from \eqref{e56fb} we can infer that
\begin{equation}\label{e57ffb}
\frac{1}{V(x)}\Delta Z(x)\geq - c K |x|^{-2 +\alpha}\geq -1\quad \forall\,\,x\in \mathbb Z^n, |x|\geq 1,
\end{equation}
provided that $K\leq \frac1{c}$.
On the other hand, we also have
\begin{equation}\label{e58fb}
\frac 1{V(0)}\Delta Z(0)
= -K\frac{\Delta(\log(2+r))(0)}{V(0)}\geq  -1,
\end{equation}
by possibly shrinking $K$.
Gathering \eqref{e57ffb} and \eqref{e58fb} we get that $Z$ satisfies \eqref{e60f}.
\end{proof}

\begin{lemma}\label{lem3z}
Let the assumptions of Theorem \ref{teo2z} hold. Then, for some $\sigma>0, K>0, \gamma>0$,
\[h(x):=\frac{\sigma}{(K+|x|^2)^{\gamma}} \quad (x\in \mathbb Z^n)\]
fulfills
\[\frac1{V(x)} \Delta h(x) \leq -1 \quad \text{for every $x\in\mathbb{Z}^n$
with $|x|\geq 1$}.\]
\end{lemma}
\begin{proof}
By arguing as in the proof of \cite[Theorem 6.1]{MoPuSo} we have that
\begin{equation}\label{e61f}
\begin{aligned}
\Delta h(x) &\leq \frac{\sigma}{(K+|x|^2)^{\gamma+2}}\Big\{-\left(\gamma-\frac{2\gamma(\gamma+1)}{n}\right)|x|^2\\
&\qquad\qquad -\left(\gamma K -\frac{\gamma(\gamma+1)}{2} \right)+c(|x|+1)\Big\} \quad
\forall\,\,x\in \mathbb Z^n,
\end{aligned}
\end{equation}
for some `absolute' constant $c > 0$.
Now we take $0<\gamma<\frac{n-2}{2}, K>\frac{\gamma+1}2$, so
$$\delta_0:=\gamma-\frac{2\gamma(\gamma+1)}{n}>0.$$
Hence, from \eqref{e61f} and \eqref{e12fb}, we can deduce that
\begin{equation}\label{e62f}
\begin{aligned}
\frac1{V(x)}\Delta h(x) &\leq - \frac{\sigma c}{(K+|x|^2)^{\gamma+2}}|x|^2(1+|x|)^{\alpha}\\
&\leq -\sigma c |x|^{\alpha-2-2\gamma}\leq -1 \quad \text{for every
$x\in \mathbb Z^n$ with $|x|\geq 1$},
\end{aligned}
\end{equation}
provided that $\gamma\leq \frac{\alpha-2}{2}$ and $\sigma>{1}/{c}\,.$
\end{proof}

\subsection{Proof of Theorems \ref{teo1z}, \ref{teo2z}}

\begin{proof}[Proof of Theorem \ref{teo1z}]
The conclusion follows by Theorem \ref{teo1}, with $d$ being the Euclidean distance, combined with Lemma \ref{lem1z} when $\alpha\in [0, 2)$, and with Lemma \ref{lem2z} when $\alpha\in [0, 2].$
\end{proof}

\begin{proof}[Proof of Theorem \ref{teo2z}]
The conclusion follows by Proposition \ref{prop71}, with $d$ being the Euclidean distance, combined with Lemma \ref{lem3z}.
\end{proof}
\bigskip
\bigskip

\noindent{\bf Acknowledgement}.
All authors are member of the ``Gruppo Nazionale per l'Analisi Ma\-te\-ma\-tica, la Probabilit\`a e le loro
Applicazioni'' (GNAMPA) of the ``Istituto Nazionale di Alta Matematica'' (INdAM, Italy), and are partially
supported by GNAMPA. The first author is partially supported by the PRIN 2022 project 2022R537CS \emph{$NO^3$ -
Nodal Optimization, NOnlinear elliptic equations, NOnlocal geometric problems, with a focus on regularity},
founded by the European Union - Next Generation EU.
The second author acknowledges that this work is part of the PRIN project 2022 Geometric-analytic methods for PDEs
and applications, ref. 2022SLTHCE, financially supported by the EU, in the framework of the "Next Generation EU
initiative".



\end{document}